\documentclass[12TP,draft]{article}

\begin{document}

\begin{center}
\LARGE\noindent\textbf{ On cyclability  of  digraphs }\\

\end{center}
\begin{center}
\noindent\textbf{Samvel Kh. Darbinyan }\\

Institute for Informatics and Automation Problems, Armenian National Academy of Sciences

E-mails: samdarbin@ipia.sci.am\\
\end{center}

\textbf{Abstract}

Given a directed graph $D$ of order $n\geq 4$  and a nonempty subset $Y$ of vertices of $D$ such that in $D$ every vertex of $Y$ reachable from every other vertex of $Y$. Assume that for every triple $x,y,z\in Y$ such that $x$ and $y$ are nonadjacent: If there is no arc from $x$ to $z$, then $d(x)+d(y)+d^+(x)+d^-(z)\geq 3n-2$. If there is no arc from $z$ to $x$, then $d(x)+d(y)+d^+(z)+d^-(x)\geq 3n-2$. We prove that there is a directed  cycle in $D$ which contains all the vertices of $Y$, except possibly one. This result is best possible in some sense and gives a answer to a question of H. Li, Flandrin and Shu (Discrete Mathematics, 307 (2007) 1291-1297).\\


\section {Introduction} 

For convenience of the  reader, terminology and notations will be given in details in section 2. A set $S$ of vertices in a directed graph $D$ (an undirected graph $G$) is said to be cyclable in $D$ (in $G$) if $D$ (if $G$) contains a directed cycle (undirected cycle) through all the vertices of $S$. There are many well-known conditions which guarantee the cyclability of a set of vertices in an undirected graph. Most of them can be seen as restrictions of Hamiltonian conditions to the considered set of vertices (See \cite{ [5], [12], [16], [17], [19]}). Let us cite for example the following:\\

 \noindent\textbf{Theorem A} ({\it R. Shi \cite{[17]}}). {\it Let $G$ be a 2-connected undirected graph of order $n$. If $S$ is a subset of the vertices of $G$ and $d(x)\geq n/2$ for all vertices $x\in S$, then $S$ is cyclable in $G$.}\\ 

 \noindent\textbf{Theorem B} ({\it R. Shi \cite{[17]}}). {\it Let $G$ be a 2-connected undirected graph of order $n$. If $S$ is a subset of the vertices of $G$ and $d(x)+d(y)\geq n$ for any two nonadjacent vertices $x\in S$ and $y\in S$, then $S$ is cyclable in $G$.}\\ 

Notice that  Theorems A and B generalize the classical theorems on hamiltonicity of Dirac and Ore, respectively.

For general directed graphs (digraphs) there are not in literature as many conditions as for undirected graphs that guarantee the existence of a directed cycle with given properties (in particular, sufficient conditions for the existence of a Hamiltonian cycles in a digraphs). 
 The more general and classical ones is the following theorem of M. Meyniel:\\ 

\noindent\textbf{Theorem C} ({\it M. Meyniel \cite{[14]}}). {\it If $D$ is a strongly connected digraph of order $n\geq 2$ and $d(x)+d(y)\geq 2n-1$ for all pairs of nonadjacent vertices $x$ and $y$ of $D$, then $D$ is Hamiltonian.}\\

Notice  that Meyniel's theorem is a common generalization of well-known  classical theorems of Ghouila-Houri \cite{[11]} and Woodall \cite{[21]}. A beautiful short proof  Meyniel's theorem can be found in \cite{[6]} (see also \cite{[20]}, pp.399-400). \\

In \cite{[8]} it was proved the following:

\noindent\textbf {Theorem D} ({\it S. Darbinyan \cite{[8]}}). {\it Let $D$ be a strongly connected digraph of order $n\geq 3$. If $d(x)+d(y)\geq 2n-1$ for any two nonadjacent vertices $x,y\in V(D)\setminus \{z_0\}$, where $z_0$ is  some vertex of $D$, then $D$ is Hamiltonian or contains a cycle of length $n-1$.}\\

The following results are  corollaries of Theorem D.\\

\noindent\textbf{Corollary 1}. {\it Let $D$ be a strongly connected digraph of order $n\geq 3$. If $D$ has $n-1$ vertices of degree at least $n$, then $D$ is  Hamiltonian or contains a cycle of length $n-1$.}\\

\noindent\textbf{Corollary 2}. {\it Let $D$ be a strongly connected digraph of order $n\geq 3$, which satisfies the conditions of Theorem D. Then $D$ has a cycle that contains all the vertices of $D$ may be except $z_0$.} \\

Let $D$ be a digraph of order $n\geq 3$.
A Meyniel set $M$ is a subset of $V(D)$ such that  $d(x)+d(y)\geq 2n-1$ for every pair of distinct vertices $x$, $y$ in $M$ which are nonadjacent in $D$. A sufficient condition for cyclability in digraphs with the condition of Meyniel's theorem was given by   K. A. Berman and X. Liu \cite{[4]}. They improved Theorem D proving the following generalization of well-known theorem of Meyniel.\\

  \noindent\textbf{Theorem E} ({\it K. Berman and X. Liu \cite{[4]}}). {\it Let $D$ be a strongly connected digraph of order $n$. Then every Meyniel set $M$ of $D$ lies in a directed cycle.}\\

Theorem E also generalizes the classical theorems of A. Ghouila-Houri \cite{[11]} and D.R. Woodall \cite{[21]}.\\

In view of the next theorem we need the following definition. Let $D$ be a directed graph and let $S$ be a nonempty subset of vertices of $D$. Following \cite{[12]}, we say that a digraph $D$ is $S$-strongly connected if for any pair $x,y$ of distinct vertices of $S$ there exists a path from $x$ to $y$ and a path from $y$ to $x$ in $D$.\\

Later H. Li, E. Flandrin and J. Shu \cite{[12]} proved the following generalization of Theorem E.\\

\noindent\textbf{Theorem F} ({\it H. Li, E. Flandrin and  J. Shu \cite{[12]}}). {\it Let $D$ be a digraph of order $n$ and $M$ be a Meyniel set in $D$. If $D$ is $M$-strongly connected, then $D$ contains a cycle through all the vertices of  $M$.} \\

Let $D$ be a digraph of order $n$. 
We say that a nonempty subset $Y$ of the vertices of $D$ satisfies  condition $A_0$ if for every triple of the vertices $x,y,z$ in $Y$ such that $x$ and $y$ are nonadjacent: if there is no arc from $x$ to $z$, then $d(x)+d(y)+d^+(x)+d^-(z)\geq 3n-2$. If there is no arc from $z$ to $x$, then $d(x)+d(y)+d^-(x)+d^+(z)\geq 3n-2$.\\
 
Y. Manoussakis \cite{[13]} proved a sufficient condition  for hamiltonicity of digraphs that involves triples rather than pairs of vertices.\\

\noindent\textbf{Theorem G} ({\it Y. Manoussakis \cite{[13]}}). {\it Let $D$ be a strongly connected digraph $D$ of order $n\geq 4$. If $V(D)$ satisfies  condition $A_0$, then $D$ is  Hamiltonian.}  \\
 
H. Li, Flandrin and  Shu \cite{[12]} (see also B. Ning \cite{[15]})  was put a question to know if this  theorem of Manoussakis (or the sufficient conditions of hamiltonicity of digraphs of Bang-Jensen, Gutin and Li \cite{[2]} or of Bang-Jensen, Guo and Yeo \cite{[3]}) has a cyclable version.\\

In this paper we prove the following theorem which gives some answer for the above question when a subset $Y\not=\emptyset$ of the vertices of a digraph $D$ satisfies  condition $A_0$ and the digraph $D$ is $Y$-strongly connected. \\

\noindent\textbf{Theorem}. {\it Let $D$ be a  digraph of order $n\geq 4$ and let $Y$ be a nonempty subset of the vertices of $D$.  Suppose that $D$ is $Y$-strongly connected and the subset $Y$ satisfies  condition $A_0$. Then $D$ contains a cycle through all the vertices  of $Y$ may be except one.}  \\

{\it Remark 1}. The following example shows that there is a digraph $D$ which contains a nonempty subset $Y$ of $V(D)$ such that $D$ is $Y$-strongly connected and the subset $Y$ satisfies  condition $A_0$ but $D$ has no cycle that contains all the vertices of $Y$.

To see this, let $G$ and $H$ be two arbitrary disjoint digraphs with $|V(G)|=m\geq 2$ and $|V(H)|=n-m\geq 4$. Let $y\in V(H)$ and $x,z\in V(G)$, $x\not=z$. Assume that $d(y,H)=2(n-m-1)$, $G$ contains a Hamiltonian cycle, $d^+(x,G)=m-1$ and $d(z,G)=2(m-1)$. From $G$ and $H$ we form a new digraph $D$ with $V(D)=V(G)\cup V(H)$ as follows: add the all possible arcs $ux, xu$, where $u\in V(H)\setminus \{y\}$, and the arc $yx$. An easy computation shows that 
$$
d(y)+d(z)+d^-(y)+d^+(x)= 4n-m-6\geq 3n-2,
$$
since $m\leq n-4$. Thus we have that the set $Y=\{x,y,z\}$ satisfies condition $A_0$, $D$ is $Y$-strongly connected and has no cycle that contains all the vertices of $Y$.\\

 Our proofs are based on the arguments of \cite{[12], [13]}.\\

\section {Terminology and Notations}

We shall assume that the reader is familiar with the standard
terminology on the directed graphs (digraph)
 and refer the reader to \cite{[1]} for terminology not discussed here.
  In this paper we consider finite digraphs without loops and multiple arcs. 
 For a digraph $D$, we denote
  by $V(D)$ the vertex set of $D$ and by  $A(D)$ the set of arcs in $D$. The order of $D$ is the number
  of its vertices. 
 The arc of a digraph $D$ directed from
   $x$ to $y$ is denoted by $xy$ or $x\rightarrow y$. If $x,y,z$ are distinct vertices in $D$, then $x\rightarrow y\rightarrow z$ denotes that $xy$ and $yz\in A(D)$. Two distinct vertices $x$ and $y$ are adjacent if $xy\in A(D)$ or $yx\in A(D) $ (or both). If $x\in V(D)$
   and $A=\{x\}$ we write $x$ instead of $\{x\}$. The out-neighborhood of a vertex $x$ is the set $N^+(x)=\{y\in V(D) / xy\in A(D)\}$ and $N^-(x)=\{y\in V(D) / yx\in A(D)\}$ is the in-neighborhood of $x$. Similarly, if $A\subseteq V(D)$, then $N^+(x,A)=\{y\in A / xy\in A(D)\}$ and $N^-(x,A)=\{y\in A / yx\in A(D)\}$.
 The out-degree of $x$ is $d^+(x)=|N^+(x)|$ and $d^-(x)=|N^-(x)|$ is the in-degree of $x$. Similarly, $d^+(x,A)=|N^+(x,A)|$ and $d^-(x,A)=|N^-(x,A)|$. 
The degree of the vertex $x$ in $D$ is defined as $d(x)=d^+(x)+d^-(x)$ (similarly, $d(x,A)=d^+(x,A)+d^-(x,A)$). 
The path (respectively, the cycle) consisting of the distinct vertices $x_1,x_2,\ldots ,x_m$ ( $m\geq 2 $) and the arcs $x_ix_{i+1}$, $i\in [1,m-1]$  (respectively, $x_ix_{i+1}$, $i\in [1,m-1]$, and $x_mx_1$), is denoted  by $x_1x_2\cdots x_m$ (respectively, $x_1x_2\cdots x_mx_1$). The length of a cycle or path is the number of its arcs.
We say that $x_1x_2\cdots x_m$ is a path from $x_1$ to $x_m$ or is an $(x_1,x_m)$-path. An $(x,y)$-path $P$ is an $(X,Y)$-path if $x\in X$, $y\in Y$ and $V(P)\cap (X\cup Y)= \{x,y\}$, where $X$ and $Y$ are subset of the  vertices of a digraph $D$.

Given a vertex $x$ of a directed path $P$ or of a directed cycle $C$, we use the notations $x^+$ and $x^-$ for the successor and the predecessor of $x$ (on $P$ or on $C$) according to the orientation, and in case of ambiguity, we precise $P$ or $C$ 
a subscript (that is $x^+_P$ ...).

A cycle (respectively, a path) that contains  all the vertices of $D$ is a  Hamiltonian cycle (respectively, is a hamiltonian path). 
A digraph is Hamiltonian if it contains a Hamiltonian cycle.
  For a cycle  $C:=x_1x_2\cdots x_kx_1$ of length $k$, the subscripts considered modulo $k$, i.e., $x_i=x_s$ for every $s$ and $i$ such that  $i\equiv s\, (\hbox {mod} \,k)$. 
If $P$ is a path containing a subpath from $x$ to $y$ we let $P[x,y]$ denote that subpath. Similarly, if $C$ is a cycle containing vertices $x$ and $y$, $C[x,y]$ denotes the subpath of $C$ from $x$ to $y$. If $C$ is a cycle and $P$ be a path in a digraph $D$, often we will write $C$ instead of $V(C)$ and $P$ instead of $V(P)$.
A digraph $D$ is strongly connected (or, just, strong) if there exists a path from $x$ to $y$ and a path from $y$ to $x$ for every pair of distinct vertices $x,y$.

Let $C$  be a non-Hamiltonian cycle in a digraph $D$. For the cycle $C$, a $C$-bypass is a path of length at least two with both end-vertices on $C$ and no other vertices on $C$. If $(x,y)$-path $P$ is a $C$-bypass with  $V(P)\cap V(C)=\{x,y\}$, 
then we call the length of the path $C[x,y]$ the gap of $P$ with respect to $C$. 

If we consider a subset of vertices $S\subseteq V(D)$, we denote the vertices of $S$ by $S$-vertices and the number of $S$-vertices in a cycle is called its $S$-length. 

The subdigraph of a digraph $D$ induced by a subset $A$ of $V(D)$ is denoted by 
$D\langle A\rangle$, or $\langle A\rangle$ for brevity.
 
 For an undirected graph $G$, we denote by $G^*$ symmetric digraph obtained from $G$ by replacing every edge $xy$ with the pair $xy$, $yx$ of arcs. 
We denote  the complete undirected graph on $n$ vertices (respectively, undirected complete bipartite graph, with partite sets of cardinalities $n$ and $m$) by $K_n$ (respectively, by $K_{n,m}$), and  $\overline K_n$ denotes the  complement of $K_n$.   If $G_1$ and $G_2$ are undirected graphs, then $G_1\cup G_2$ is the disjoint union of $G_1$ and $G_2$. The join of $G_1$ and $G_2$, denoted by  $G_1 + G_2$, is the  union of $G_1\cup G_2$ and of all the edges between $G_1$ and $G_2$. The converse digraph $\overleftarrow {D}$  of a digraph $D$ is the digraph obtained from $D$ by reversing all arcs of $D$.

 For integers $a$ and $b$, $a\leq b$, let $[a,b]$  denote the set of
all integers which are not less than $a$ and are not greater than
$b$.  

\section { Preliminaries }

We now collect the tools which we need in  proof of our theorem.
In the following, we often use the following definition: 

\noindent\textbf{Definition}. {\it Let $P=x_1x_2\ldots x_m$ ($m\geq 2$) be a path in a digraph $D$ and let $Q=y_1y_2\ldots y_k$ be a path in $\langle V(D)\setminus V(P)\rangle$ (possibly, $k=1$). Assume that there is an $i\in [1,m-1]$ such that $x_iy_1$ and $y_kx_{i+1}\in A(D)$. In this case $D$ contains a path $x_1x_2\ldots x_iy_1y_2\ldots y_k x_{i+1}\ldots x_m$ and we say that $Q$ can be inserted into $P$.}\\

The following Lemmas 3.1 and 3.2 are slight modifications of lemma by H\"{a}ggkvist and Thomassen \cite{[10]} and of lemma by Bondy and Thomassen \cite{[6]}, respectively (their proofs  are not too difficult).
 They
will be used extensively in the proof of our result. \\

\noindent\textbf{Lemma 3.1}. {\it Let $C_k:=x_1x_2\ldots x_kx_1$, $k\geq 2$, be a non-Hamiltonian cycle in a digraph $D$. Moreover, assume that there exists a path $Q:=y_1y_2\ldots y_r$, $r\geq 1$, in $\langle V(D)\setminus V(C_k)\rangle$. If 
 $d^-(y_1,C_k)+d^+(y_r,C_k)\geq k+1$,
 then for all $m\in [r+1,k+r]$ the digraph $D$ contains a cycle $C_m$ of length $m$  with vertex set $V(C_m)\subseteq V(C_k)\cup V(Q)$}.  \fbox \\\\

\noindent\textbf{Lemma 3.2}. {\it Let $P:=x_1x_2\ldots x_k$, $k\geq 2$, be a non-Hamiltonian path in a digraph $D$. Moreover, assume that there exists a path $Q:=y_1y_2\ldots y_r$, $r\geq 1$, in $\langle V(D)\setminus V(P)\rangle$. If

 $$d^-(y_1,P)+d^+(y_r,P)\geq k+d^-(y_1,\{x_k\})+d^+(y_r,\{x_1\}),$$ 
then there is an $i\in [1,k-1]$ such that $x_iy_1$ and $y_rx_{i+1}\in A(D)$, i.e., $D$ contains a path from $x_1$ to $x_k$ with vertex set $V(P)\cup V(Q)$, i.e., $Q$ can be inserted into $P$}. \fbox \\\\

The following lemma from \cite{[12]} is a slight  modification of Multi-Insertion Lemma due to Bang-Jensen, Gutin and H. Li  (see \cite{[1]}, Lemma 5.6.20).

\noindent\textbf{Lemma 3.3} ({\it H. Li, E. Flandrin. J. Shu  \cite{[12]}}). {\it  Let $D$ be a digraph and let $P$ be an $(a,b)$-path in $D$. Let $Q$ be a 
 path in $\langle V(D)\setminus V(P)\rangle$ and let $S$ be a subset of $V(Q)$. If every vertex of $S$ can be inserted into $P$, then there exists an $(a,b)$-path $R$ such that $V(P)\cup S \subseteq V(R)\subseteq V(P)\cup V(Q)$.} \fbox \\\\

The following lemma also  was proved in \cite{[12]}.

\noindent\textbf{Lemma 3.4} ({\it H. Li, E. Flandrin. J. Shu \cite{[12]}}). {\it Let $D$ be a digraph of order $n$ and $S\subset V(D)$, $S\not= \emptyset$. Assume that $D$ is $S$-strongly connected and satisfies for any pair of nonadjacent vertices $x$, $y$ in $S$ the degree condition $d(x)+d(y)\geq 2n-1$. If $C$ is a cycle in $D$ of maximum $S$-length and $s$ is an $S$-vertex of $V(D)\setminus V(C)$, then $D$ contains a $C$-bypass through $s$.} \fbox \\\\

By inspection of the proof in \cite{[12]} one can state Lemma 3.4 in the following form (its proof is the same as the proof of Lemma 3.4).\\ 

\noindent\textbf{Lemma 3.5}. {\it Let $D$ be a digraph of order $n$ and let $C$ be a non-Hamiltonian cycle in $D$. Let $x$ be an arbitrary  vertex not on $C$. Assume that in $D$ there are $(C,x)$- and $(x,C)$-paths and $D$ contains no $C$-bypass through $x$. Then the following holds:

(i). If $x$ is adjacent to some vertex $y$ of $C$, then $D$ is not 2-strong, $d(x,V(C)\setminus \{y\})=0$ and $d(x)+d(z)\leq 2n-2$ for all the vertices $z\in V(C)\setminus \{y\}$.

(ii). Assume that $x$ and any vertex of $C$ are nonadjacent, i.e., $d(x,V(C))=0$. Let $P$ be a shortest $(C,x)$-path with $\{u\}= V(P)\cap V(C)$ and let $Q$ be a shortest $(x,C)$-path with $\{v\}= V(Q)\cap V(C)$. 
Then the following holds: 

If $u\not= v$, then $d(x)+d(z)\leq 2n-2$ for all the vertices $z\in V(C)$, may be except one from $\{u,v\}$.
 
If $u=v$, then $d(x)+d(z)\leq 2n-2$ for all the vertices $z\in V(C)\setminus \{u\}$}. \fbox \\\\

In \cite{[13]} it was proved the following

\noindent\textbf{Lemma 3.6} ({\it Y. Manoussakis \cite{[13]}}). {\it Let  $D$ be a digraph of order $n$ and let $V(D)$ satisfies  condition $A_0$. Assume that there are two distinct pairs of nonadjacent vertices $x,y$ and $x,z$ in $D$. Then either $d(x)+d(y)\geq 2n-1$ or $d(x)+d(z)\geq 2n-1$.} \fbox \\\\

It is not difficult to show that one can Lemma 3.6 state the following much stronger form:\\

\noindent\textbf{Lemma 3.7.} {\it Let  $Y$ be a subset of vertices in a digraph $D$ of order $n$ and let $Y$  satisfies  condition $A_0$. Assume that there are two distinct pairs of nonadjacent vertices $x,y$ and $x,z$ in $Y$. Then either $d(x)+d(y)\geq 2n-1$ or $d(x)+d(z)\geq 2n-1$.} \fbox \\\\

For the proof of our results we also need the following simple lemma.

\noindent\textbf{Lemma 3.8}. {\it Let  $D$ be a  digraph of order $n$. Assume that $xy\notin A(D)$ and the vertices $x$, $y$ in $D$  satisfies the degree condition  $d^+(x)+d^-(y)\geq n-2+k$, where $k\geq 1$. Then $D$ contains at least $k$ internally disjoint $(x,y)$-paths of length two.} \fbox \\\\

The following lemma was also proved in \cite{[12]}.

\noindent\textbf{Lemma 3.9} ({\it H. Li, E. Flandrin. J. Shu  \cite{[12]}}). {\it Let  $D$ be a  digraph of order $n$ and $S\subseteq V(D)$, $S\not=\emptyset$. Assume that $D$ is $S$-strongly connected and satisfies for any  pair of nonadjacent vertices $x,y$ in $S$  the degree condition  $d(x)+d(y)\geq 2n-1$. Then any two $S$-vertices $s$ and $s'$ are contained in a cycle of $D$ such that they are at distance at most two on this cycle.} \fbox \\\\

One can Lemma 3.9 state the following form.\\

\noindent\textbf{Lemma 3.10} . {\it Let  $D$ be a  digraph of order $n$. Assume that a pair of distinct vertices $x,y$ in $D$ satisfies the degree condition  $d(x)+d(y)\geq 2n-1$.
If $D$ is $\{x,y\}$-strongly connected, then the vertices $x$ and $y$ are contained in a cycle of $D$ such that they are at distance at most two on this cycle.} \fbox \\\\

Now we will prove the following lemma.\\

\noindent\textbf{Lemma 3.11}. {\it Let  $D$ be a digraph of order $n$ and let $Y$ be a subset of vertices of $D$ with $|Y|\geq 4$. Assume that $D$ is $Y$-strongly connected and the subset $Y$ satisfies  condition $A_0$. If $C$ is a non-Hamiltonian cycle in $D$ which contains at least two $Y$-vertices, then for every  $Y$-vertex $y$  of $V(D)\setminus V(C)$ there is a $C$-bypass through $y$.}

\noindent\textbf{Proof of Lemma 3.11}. If the cycle $C$ contains at least three $Y$-vertices, then the lemma  immediately follows from Lemmas 3.5 and 3.7. Assume therefore that $C$ contains exactly two $Y$-vertices, say $x$ and $u$, and there exists a $Y$-vertex, say $y$, in 
$B:=V(D)\setminus V(C)$ such that in $D$ there is no $C$-bypass through $y$. From $|Y|\geq 4$ it follows that $B$ contains at least two $Y$-vertices. Let $z$ be an arbitrary $Y$-vertex of $B$ other than $y$.

 We will consider two cases.\\

{\it Case 1}. $d(y,C)\geq 1$.

Without loss of generality, assume that the vertex $y$ is adjacent to a vertex $w$ of  $ V(C)$. If $w\notin \{u,x\}$, then from Lemma 3.5(i) it follows that $y,u$ and $y,x$ are distinct pairs of nonadjacent vertices of $Y$,  $d(y)+d(u)\leq 2n-2$ and $d(y)+d(x)\leq 2n-2$, which contradicts Lemma 3.7. Assume therefore that $w\in \{x,u\}$, for example, let $w=u$ and $yu\in A(D)$. Since we assumed that $D$ has no $C$-bypass through $y$, by Lemma 3.5(i) we have
$$
d(y)+d(x)\leq 2n-2 \quad \hbox{and} \quad d(y,V(C)\setminus \{u\})=0. \eqno (1)
$$

Now we distinguish two subcases.\\

{\it Subcase 1.1}. $xz\in A(D)$.

Then $zy\notin A(D)$ (for otherwise $xzyu$ is a $C$-bypass through $y$, which contradicts the our assumption that $D$ has no $C$-bypass through $y$). Therefore, the  triple of $Y$-vertices $x, y,z$ satisfies condition $A_0$, i.e.,
$$
d(y)+d(x)+d^-(y)+d^+(z)\geq 3n-2.
$$
This together with $d(y)+d(x)\leq 2n-2$ (by (1)) implies that $d^+(z)+d^-(y)\geq n$. Hence, by Lemma 3.8, $z\rightarrow v\rightarrow y$ for some vertex $v$ other than $u$. 
From $d(y,V(C)\setminus \{u\})=0$ (by (1)) it follows that $v\in B$. Thus, $xzvyu$ is a $C$-bypass through $y$, a contradiction.\\

{\it Subcase 1.2}. $xz\notin A(D)$.

Then by condition $A_0$ we have
$$
d(y)+d(x)+d^+(x)+d^-(z)\geq 3n-2.
$$
Therfore, by (1), $d^+(x)+d^-(z)\geq n$, and hence by Lemma 3.8 and $xy\notin A(D)$, there exists a vertex $v$ other that $u$ and $y$ such that 
$x\rightarrow v\rightarrow z$. It is easy to see that $vy\notin A(D)$.
 Again we have, $zy\notin A(D)$ and $d^+(z)+d^-(y)\geq n$. Hence, by Lemma 3.8, $z\rightarrow a\rightarrow y$ for some vertex $a$ other than $u$. It is not difficult to see that $a\in B\setminus \{y,z,v\}$. Consequently, $vzayu$ or $xvzayu$ is a $C$-bypass through $y$ when $v\in C$ or not, respectively, a contradiction. The discussion of Case 1 is completed.\\

{\it Case 2}. $d(y,C)=0$.

By Lemma 3.5(ii), we have either 
$d(y)+d(x)\leq 2n-2$ or $d(y)+d(u)\leq 2n-2$. Without loss of generality, assume that 
$$
d(y)+d(x)\leq 2n-2. \eqno (2)
$$
This together with condition $A_0$  implies that 
$$
d^+(y)+d^-(u)\geq n \quad \hbox{and} \quad d^+(u)+d^-(y)\geq n. \eqno (3)
$$
 This together with Lemma 3.8 implies that there are vertices $a$ and $v$ (possibly, $a=v$) other than $z$ such that $u\rightarrow v\rightarrow y$ and $y\rightarrow a\rightarrow u$. Observe that $v$ and $a$ are not on $C$ since $d(y,C)=0$. 

Assume first that there is a vertex $w\in V(C)\setminus \{u\}$ which is adjacent to $z$. Without loss of generality, assume that $zw\in A(D)$ (for the case $wz\in A(D)$ we will consider the converse digraph of $D$). 
If $yz\in A(D)$, then $uvyzw$ is a $C$-bypass through $y$, a contradiction. Assume therefore that $yz\notin A(D)$. Then from (2) and  condition $A_0$ it follows that $d^+(y)+d^-(z)\geq n$. 
Therefore, by Lemma 3.8, for some vertex $b\in B\setminus \{v\}$,  $y\rightarrow b\rightarrow z$, and hence, $uvybzw$ is a $C$-bypass through $y$, which is a contradiction.

Assume second that $d(z,V(C)\setminus \{u\})=0$. In particular, the vertices $z$ and $x$ are nonadjacent. From (2) and  condition $A_0$ it follows that 
$$
d^+(x)+d^-(z)\geq n \quad \hbox{and} \quad d^+(z)+d^-(x)\geq n.
$$
From $d^+(z)+d^-(x)\geq n$ and  Lemma 3.8 it follows that there are at least two $(z,x)$-paths of length two.\\

 {\it Subcase 2.1}. There is a $(z,x)$-paths of length two, say $z\rightarrow b\rightarrow x$, such that $b\notin \{u,v\}$.
 
 Then $yz\notin A(D)$ and $yb\notin A(D)$ (for otherwise, $uvyzbx$ or $uvybx$ is a $C$-bypass through $y$, when $yz\in A(D)$ and $yb\in A(D)$, respectively). Since $x,y,z$ are $Y$-vertices, from condition $A_0$ and (2) it follows that $d^+(y)+d^-(z)\geq n$.
 Now using Lemma 3.8 and the facts that $yz\notin A(D)$ and $yb\notin A(D)$, we obtain that there exists a vertex $q\in B\setminus \{v,b,z\}$ such that $y\rightarrow q\rightarrow z$. Thus, $uvyqzbx$ is a $C$-bypass through $y$, a contradiction.\\

{\it Subcase 2.2}. There is no $w\in B\setminus \{v\}$ such that $z\rightarrow w\rightarrow x$. 

Then from Lemma 3.8 and $d^+(z)+d^-(x)\geq n$ it follows that $d^+(z)+d^-(x)= n$  and $zv,vx,zu\in A(D)$ (i.e., there are exactly two $(z,x)$-paths of length two). 
Now using the inequality $d^+(u)+d^-(y)\geq n$ (by (3)) and Lemma 3.8 we conclude that there exist at least two $(u,y)$-paths of length two. If there is a path $u\rightarrow c\rightarrow y$    such that $c$ is other than $v$ and $z$,  then we may consider the paths  $u\rightarrow c\rightarrow y$ and $z\rightarrow v\rightarrow x$.
 For these paths we have the above considered case ($b\notin \{u,v\}$). 
Assume therefore that there is no 
$c\in B\setminus \{v,z\}$ such that $u\rightarrow c\rightarrow y$. Again using Lemma 3.8, it is easy to see that  $u\rightarrow z\rightarrow y$. From this and Lemma 3.8 it follows that $d^+(u)+d^-(y)= n$ since $d^+(u)+d^-(y)\geq n$ (by (3)). Now by  condition $A_0$ and (2) we have
 $$
d(y)+d(x)= 2n-2.
$$

If $yv\in A(D)$, then $uzyvx$ is a $C$-bypass through $y$, a contradiction. Assume therefore that $yv\notin A(D)$. This together with $d^+(y)+d^-(u)\geq n$ (by (3)) implies that there is a vertex $u_1\in B\setminus \{z,v\}$ such that $y\rightarrow u_1\rightarrow u$.  
From this it is easy to see that $xv\notin A(D)$ (for otherwise $xvyu_1u$ is a $C$-bypass through $y$). 
This together with $d^+(x,\{y,z,v\})=0$ and $d^+(x)+d^-(z)\geq n$ implies that $x\rightarrow u_1\rightarrow z$ (for otherwise there is a vertex $r\in B\setminus \{v,y,z,u_1\}$ such that $x\rightarrow r \rightarrow z$ and hence, $xrzyu_1u$ is a $C$-bypass through $y$, a contradiction).

 If $d^+(x)+d^-(y)\geq n$, then there is a vertex $u_2\in B\setminus \{v,z, u_1\}$ such that $x\rightarrow u_2\rightarrow y$ since $d^+(x,\{v,y,z\})=0$ and $uy\notin A(D)$. 
Therefore, $xu_2yu_1u$ is a $C$-bypass through $y$.
 Assume therefore that  $d^+(x)+d^-(y)\leq n-1$. 
Then from $d(x)+d(y)=2n-2$ it follows that 
 $d^+(y)+d^-(x)\geq n-1$. Now using Lemma 3.8 and the  facts that $d^+(y, V(C)\cup \{v\})=0$ and $zx\notin A(D)$ we obtain that $y\rightarrow w\rightarrow x$ for some $w\in B\setminus \{v,z\}$. Therefore $uvywx$ is a $C$-bypass through $y$, which is a contradiction.
 Lemma 3.11 is proved. \fbox \\\\

\section {Proof of the main result}

For readers convenience, again we will formulate the main result.

\noindent\textbf{Theorem}. {\it Let $D$ be a  digraph of order n and let $Y$ be a nonempty subset of the vertices of $D$, where $|Y|\geq 2$.  Suppose that $D$ is $Y$-strongly connected and the subset $Y$ satisfies  condition $A_0$. Then $D$ contains a cycle through all the vertices  of $Y$ may be except one.} 

\noindent\textbf{Proof of the theorem }. Suppose, on the contrary, that is a digraph $D$ and a nonempty subset $Y$ of the vertices of $D$ satisfy the supposition of the theorem  but any cycle in $D$ does not contain at least two $Y$-vertices. By  Manoussakis' theorem, we may assume that $Y\not= V(D)$.
Since $D$ is $Y$-strongly connected, using Lemmas 3.7 and 3.10 we obtain that $|Y|\geq 4$ and in $D$ there exists a cycle which contains at least two $Y$-vertices. 
From Lemma 3.11 it follows that if $C$ is a non-Hamiltonian cycle in $D$, then $D$ has a $C$-bypass through every $Y$-vertex of $V(D)\setminus V(C)$.
 In $D$ we choose a cycle $C$ and a $C$-bypass $P_0$ through a $Y$-vertex of $V(D)\setminus V(C)$ such that 

(a) $C$ contains as many vertices of $Y$ as possible  ($C$ contains at least two $Y$-vertices),

(b) the gap of  $C$-bypass $P_0$ is minimum, subject to (a) (by Lemma 3.11, for the cycle $C$  there exists $C$-bypass through any $Y$-vertex not on $C$), and

(c) the length of $C$-bypass $P_0$ is minimum, subject to (a) and (b).\\ 

In the sequel we assume that the cycle $C:=x_1x_2\ldots x_mx_1$ and the $C$-bypass $P_0:=x_1z_1\ldots z_kyz_{k+1}\ldots$ $z_tx_{a+1}$ satisfy the conditions (a)-(c), where $1\leq a\leq m-1$ and $y\in V(D)\setminus V(C)$ is a $Y$-vertex 
(possibly, $z_1=y$ or $y=z_t$). 
Since the cycle $C$ has the maximum $Y$-length, it follows that $C[x_2,x_a]$ contains a $Y$-vertex. In particular, $a\geq 2$.
 Note that the gap of $C$-bypass $P_0$ is equal to $a$. 

Denote $P:=z_1\ldots z_kyz_{k+1}\ldots z_t$, $E:=P[z_1,z_k]$ and $L:=P[z_{k+1},z_t]$.  Since the gap $a$ is minimal, the vertex $y$ is not adjacent with any vertex of  $C[x_2,x_{a}]$, i.e, $d(y,C[x_2,x_a])=0$.
 Therefore, by Lemma 3.2, 
$$
d(y,C)= d(y,C[x_{a+1},x_1])\leq m-a +d^-(y,\{x_1\})+ d^+(y,\{x_{a+1}\}), \eqno (4) 
$$
since any $Y$-vertex of $B:=V(D)\setminus V(C)$ cannot be inserted into $C$.
From the minimality of the path $P$ it follows that
$$
d(y,V(P))\leq |V(P)|+1-d^-(y,\{x_1\})- d^+(y,\{x_{a+1}\}).      \eqno (5)
$$ 
Notice that $|C[x_2,x_a]|=a-1$ and $|C[x_{a+1},x_1]|=m-a+1$.
First for the cycle $C$ and $C$-bypass $P_0:=x_1z_1\ldots z_kyz_{k+1}\ldots z_tx_{a+1}$ it is convenient to prove the following two claims below.\\

\noindent\textbf{Claim 1}. {\it  There is a $Y$-vertex, say $y_1$, in $C[x_2,x_a]$ such that $y, y_1$ are nonadjacent, $d(y)+d(y_1)\leq 2n-2$ and $y_1$ cannot be inserted into $C[x_{a+1},x_1]$. Moreover,

(i). The path $P$ contains exactly one $Y$-vertex, namely only $y$.

(ii). Any $Y$-vertex of $C[x_2,x_a]$ other than $y_1$ can be inserted into $C[x_{a+1},x_1]$.

(iii). There are three $(x_{a+1},x_1)$-paths, say $P_1, P_2$ and $P_3$, with vertex set $C[x_{a+1},x_1]\cup F_1$, $C[x_{a+1},x_1]\cup F_2$ and $C[x_{a+1},x_1]\cup F_3$, respectively, where  $F_1\subseteq C[x_2,x_a]$, $F_2\subseteq C[x_2,y_1^-]$, $F_3\subseteq C[y_1^+, x_{a}]$ $($if $y_1=x_2$, then $C[x_2,y^-_1]=\emptyset$, if $y_1=x_a$, then $C[y^+_1,x_a]=\emptyset$$)$ and $F_1$ $($respectively, $F_2, F_3$ $)$ contains all the $Y$-vertices of $C[x_2,x_a]\setminus \{y_1\}$ $($respectively, all the $Y$-vertices of $C[x_2,y_1^-]$, all the $Y$-vertices of $C[y_1^+,x_a])$.}

\noindent\textbf{Proof of Claim 1}. Since we assumed that the cycle $C$ has maximum $Y$-length, from Lemma 3.3 it follows that some $Y$-vertex, say $y_1$, of $C[x_2,x_a]$ cannot be inserted into $C[x_{a+1},x_1]$. 
Hence, using Lemma 3.2, we obtain 
$$
d(y_1,C)=d(y_1,C[x_2,x_a])+d(y_1,C[x_{a+1},x_1])\leq 2a-4+m-a+2 =m+a-2.  \eqno (6)
$$

From the minimality of $C[x_2,x_a]$ it follows that the vertices $y$ and $y_1$ are nonadjacent. 

Put
$R:=V(D)\setminus (V(C)\cup V(P))$.
 Now we want to compute the sum of degree $y$ and $y_1$. By  minimality of $C[x_2,x_a]$ we have
$$
d(y_1,R)+d(y,R)\leq 2|R| \quad \hbox{and} \quad d(y,C[x_2,x_a])=0. \eqno (7)
$$
From the minimality of $C[x_2,x_a]$ also it follows that 
$$
d^+(y_1,\{z_1,\ldots , z_k\})= d^-(y_1,\{z_{k+1},\ldots , z_t\})=0. \eqno (8)
$$
Therefore, $d(y_1,P)\leq |P|-1$ since $y$ and $y_1$ are nonadjacent.
This together with the above inequalities (4)-(7) gives
$$
d(y)+d(y_1)=d(y,R)+d(y_1,R)+d(y,P)+d(y,C)+d(y_1,P)+d(y_1,C)\leq $$ 
$$ 2|R|+2|P|+2m-2=2n-2.
$$
Thus, $d(y)+d(y_1)\leq 2n-2$ for any $Y$-vertex $y$ of $P$ and for any $Y$-vertex $y_1$ of $C[x_2,x_a]$ which cannot be inserted into $C[x_{a+1},x_1]$. This together with Lemma 3.7 implies that $P$ contains only one $Y$-vertex, namely $y$, and any $Y$-vertex of $C[x_2,x_a]$ different from $y_1$ can be inserted into $C[x_{a+1},x_1]$. From this and Lemma 3.3 immediately follows the third assertion of the claim.         
Claim 1 is proved.   \fbox \\\\

\noindent\textbf{Claim 2}. {\it Let $y_1$ be a  $Y$-vertex of $C[x_2,x_a]$ which cannot be inserted  into $C[x_{a+1},x_1]$. Then $d(y_1,P)=0$.}

\noindent\textbf{Proof of Claim 2}.  Suppose, on the contrary, that $d(y_1,P)\geq 1$. Then from (8) it follows that either $z_iy_1\in A(D)$ or $y_1z_j\in A(D)$, for some $i\in [1,k]$ or $j\in [k+1,t]$, respectively. Let $y_1z_j\in A(D)$. We consider the cycle $C_1:=P_3C[x_1,y_1]P[z_j,z_t]x_{a+1}$. 
This cycle contains all the $Y$-vertices of $C$, and hence has  maximum $Y$-length. 
It is easy to see that $x_1z_1\ldots z_kyz_{k+1}\ldots z_j$ is a $C_1$-bypass through $y$. 
By choice of the cycle $C$ and $C$-bypass $P_0$ we have that  $y_1=x_a$, i.e., the $C$-gap of $P_0$ and $C_1$-gap of $Q$ are equal but the path $z_1 \ldots z_kyz_{k+1}\ldots z_{j-1}$ is short than the path $z_1 \ldots z_kyz_{k+1}\ldots z_{t}$, which contradicts (c). Therefore, $y_1z_j\notin A(D)$ for all $j\in [k+1,t]$. 
 By similar arguments, one can show that $z_iy_1\notin A(D)$ for all $i\in [1,k]$. Thus, 
$$
d^-(y_1,\{z_1,\ldots , z_k\})= d^+(y_1,\{z_{k+1},\ldots , z_t\})=0
$$
which together with (8) implies that
 $d(y_1,P)=0$. Claim 2 is proved.   \fbox \\\\

Let $x$ be an arbitrary $Y$-vertex in $B=V(D)\setminus V(C)$ other than $y$. Claim 1 implies that $x$ is not on $P$. We distinguish two cases according as in  $\langle B\setminus (P\setminus \{y\})\rangle$ there exists a path with end-vertices $x$ and $y$ or not.\\

{\it Case 1. In  $\langle B\setminus (P\setminus \{y\})\rangle$ there exists a path from $x$ to $y$ or there exists a path from $y$ to $x$.} 

Without loss of generality, we may assume that in  $\langle B\setminus (P\setminus \{y\})\rangle$ there is an $(x,y)$-path (for otherwise, we consider the converse digraph of $D$). 

Let $H$ be a shortest $(x,y)$-path in $\langle B\setminus (P\setminus \{y\})\rangle $. Observe that $x_1x\notin A(D)$, since otherwise, if $x_1x\in A(D)$, then the path $P_1$ (Claim 1(iii)) together with the arc $x_1x$ and the paths $H$ and $P_0[y,x_{a+1}]$ forms a cycle, say $C_1$, which contains more $Y$-vertices  than $C$, which contradicts to our assumption that $C$ has maximum $Y$-length ($C_1$ contains all $Y$-vertices of $C$, except $y_1$, and $Y$-vertices $x$ and $y$). 
Put $R:=B\setminus (V(P)\cup V(H))$ and $H':=H[x^+_H,y^-_H]$ (if $x^+_H=y$, then $H'=\emptyset$).
 
By Claim 1,  we have 
$$
d(y_1)+d(y)\leq 2n-2.       \eqno (9)
$$
From the minimality of the gap $a$ (or of the existence of the path $P_3$) it follows that $y_1x\notin A(D)$ (therefore either $xy_1\in A(D)$ or $x$ and $y_1$ are nonadjacent)
and 
$$
d^+(y_1,R)+d^-(x,R)\leq |R|.    \eqno (10)
$$

{\it Subcase 1.1}. $xy_1\in A(D)$. 

From Lemma 3.2 it follows that 
$$
d^-(x,P_1)+d^+(y_1,P_1)\leq |P_1|.    \eqno (11)
$$
since $x_1x\notin A(D)$ and the arc $xy_1$ cannot be inserted into $P_1$ (for otherwise, $D$ contains a cycle which contains all $Y$-vertices of $C$ and $Y$-vertices $x, y$, which is a contradiction).
 From the minimality of the gap $a$,  the existence of the paths $P_2, P_3$ (by Claim 1) and Claim 2 it follows that 
$$
d^+(y_1,P\cup H)=d^-(x,C[x_2,x_a])=d^-(x,P)=0. \eqno (12)
$$
Clearly,
$$
d^+(y_1,Q)\leq |Q|-1 \quad \hbox{and} \quad d^-(x,H')\leq |H'|,  \eqno (13)
$$
where $Q:=C[x_2,x_a]-P_1$.  By adding the above relations (10)-(13), we obtain
$$
d^+(y_1)+d^-(x)=d^+(y_1,R)+d^-(x,R)+d^-(x,P_1)+d^+(y_1,P_1)+d^+(y_1,Q)+d^+(y_1,P\cup H)+$$ $$ d^-(x,Q)+d^-(x,H')+d^-(x,P)\leq
|R|+|P_1|+|Q|+|H'|-1\leq n-2.
$$
This together with (9) gives 
$$ 
d(y)+d(y_1)+d^-(x)+d^+(y_1)\leq 3n-4,
$$
which contradicts condition $A_0$, since $y, y_1$ and  $x $ are $Y$-vertices,  $y, y_1$ are nonadjacent and $y_1x\notin A(D)$.\\

{\it Subcase 1.2. The vertices $x$ and $y_1$ are nonadjacent.}
\\
We will distinguish two subcases, according as there exists a $(y,x)$-path in $\langle B\setminus (P\setminus \{y\})\rangle$ or not.

{\it Subcase 1.2.1. In  $\langle B\setminus (P\setminus \{y\})\rangle$ there is no $(y,x)$-path, in particular $yx\notin A(D)$.}

Then, clearly
$$
d^+(y,R\cup H')+d^-(x,R\cup H')\leq |R\cup H'|=|R|+|H'|.  \eqno (14)
$$

It is not difficult to see that the path $H$ cannot be inserted into $C$. Hence, from Lemma 3.1 it follows that
$$
d^-(x,C)+d^+(y,C)\leq |C|.        \eqno (15)
$$
Moreover,
$$
d^-(x,P)=d^-(x,E)+d^-(x,L)\leq |L|, \quad \hbox{since} \quad d^-(x,E)=0,
$$
recall that  $E:=P[z_1,z_k]$ and $L:=P[z_{k+1},z_t]$, and by minimality of $P$,
$$
d^+(y,P)=d^+(y,E)+d^+(y,L)\leq |E|+1, \quad \hbox{since} \quad d^+(y,L)\leq 1.
$$
Hence
$$
d^-(x,P)+d^+(y,P)\leq |L|+|E|+1=|P|.
$$
The last inequality together with (14), (15) and (9) implies that
$$
d(y)+d(y_1)+d^-(x)+d^+(y)\leq 2n-2+|R|+|H'|+|C|+|P|=3n-3,
$$
which contradicts condition $A_0$, since $y,y_1$ and $x$ are $Y$-vertices, $y, y_1$ are nonadjacent and $yx\notin A(D)$.\\

{\it Subcase 1.2.2. In  $\langle B\setminus (P\setminus \{y\})\rangle$ there is a $(y,x)$-path.}

Assume first that in $\langle B\setminus(P\cup H\setminus \{x,y\})\rangle$ there is a $(y,x)$-path . Let $Q$ be a shortest $(y,x)$-path in $\langle B\setminus (P\cup H\setminus \{x,y\})\rangle$.

Let now $R:=B\setminus (P\cup H\cup Q)$. We want to compute the degree sum of the vertices $x$ and $y_1$. 
From the minimality of the  $C[x_2,x_{a}]$ and the existence of the paths $H$ and $Q$ it follows that
$$
d(x,R)+d(y_1,R)\leq 2|R| \quad \hbox{and} \quad d(y_1, H'\cup Q')=0, \eqno (16)
$$
where $Q':=Q[y^+_Q,x^-_Q]$ (here if $y^+_Q=x$, then $Q'=\emptyset$).
By Claim 2, $d(y_1,P)=0$. This together with (6) implies that 
$$
d(y_1,C\cup P)\leq m+a-2.    \eqno   (17)
$$

Now we consider the vertex $x$. It is not difficult to see that $x$ cannot be inserted into $P$ (for otherwise there exists an $(z_1,z_t)$-path with vertex set $V(P)\cup \{x\}$ which together with the arcs $x_1z_1, z_tx_{a+1}$ and the path  $P_1$ (Claim 1) forms a cycle which  contains all the $Y$-vertices of $C$ except $y_1$ and $Y$-vertices $x$ and $y$, this contradicts the  assumption that $C$ has the maximum $Y$-length). Therefore, by Lemma 3.2,
$$
d(x,P)\leq |P|+1.    \eqno (18)
$$
 From the minimality of the paths $H$ and $Q$ it follows that
$$
d(x,Q')\leq |Q'|+1 \quad \hbox{and} \quad d(x,H')\leq |H'|+1.  \eqno (19)
$$
Since  the gap $a$ is minimal,  we obtain that $d(x,C[x_2,x_a])=0$. Using the path $P_1$ (Claim 1), it is not difficult to see that $x_1x\notin A(D)$ and $xx_{a+1}\notin A(D)$. 
Therefore, by Lemma 3.2, $d(x,C)=d(x,C[x_{a+1},x_1])\leq m-a$, since $x$ cannot be inserted into $C$. Summing the above inequalities (16)-(19) and the last inequality, an easy computation shows that 
$$
d(y_1)+d(x)\leq 2|R|+2m+|P|+|Q'|+|H'|+1\leq 2n-2.
$$
This together with $d(y_1)+d(y)\leq 2n-2$ (by (9)) contradicts Lemma 3.7, since $y_1, x$ and $y_1, y$ are two distinct pairs of nonadjacent vertices in $Y$.

Assume second that any $(y,x)$-path in $\langle B\setminus (P\setminus\{y\})\rangle$ has a common internal vertex with $(x,y)$-path $H$. Then, in particular, the vertices $y$ and $x$ are nonadjacent. 
Let  $Q$ again be a shortest $(y,x)$-path in $\langle B\setminus (P\setminus \{y\})\rangle$. 

Denote  $Q':=Q[y^+_Q, x^-_Q]$ and $R:=B\setminus (P\cup H'\cup Q'\cup \{x\})$.
 Observe that $|H'|\geq 1$ and $|Q'|\geq 1$ since $y$ and $x$ are nonadjacent.

Now we want to compute the sum $d^+(x)+d^-(y)$.
It is easy to see that
$$
d^+(x,R)+d^-(y,R)\leq |R|,  \eqno (20)
$$
since in $\langle B\setminus (P\setminus \{y\})\rangle$ any $(y,x)$-path has a common internal vertex with the minimal path $H$ and the length of $H$ is more than or equal to two. Observe that $(y,x)$-path $Q$ cannot be inserted into $C[x_{a+1},x_1]$ (for otherwise in $D$ there is a cycle which contains more $Y$-vertices than the cycle $C$). 
Notice that $xx_{a+1}\notin A(D)$ (for otherwise, if  $xx_{a+1}\in A(D)$ then the paths $P[z_1,y]$, $Q$ and $P_1$ (Claim 1)  and the arc $xx_{a+1}$ form a cycle which has more $Y$-lengths than the cycle $C$). 
Therefore,  by Lemma 3.2 we have 
$$
d^+(x, C[x_{a+1},x_1])+d^-(y, C[x_{a+1},x_1])\leq m-a+d^-(y,\{x_{1}\})+d^+(x,\{x_{a+1}\})\leq m-a+1.  \eqno (21)
$$
From  minimality of $C[x_2,x_{a}]$ and the existence of the path $Q$ it follows that 
$$
d^-(y, C[x_{2},x_a])=d^+(x, C[x_{2},x_a])=0. \eqno (22)
$$
 By minimality of $P$ we have 
$$
d^-(y,P)=d^-(y,E)+d^-(y,L)\leq |L|+1. 
$$
 On the other hand
$$
 d^+(x,P)=d^+(x,E)+d^+(x,L)\leq |E|,  
$$
since if $xz_j\in A(D)$ for some $j\in [k+1,t]$, then using the paths $Q$, $P_1$ and the subpaths of the path $P$ we can obtain  a cycle which contains more $Y$-vertices than $C$. The last two inequalities imply that
$$
d^-(y,P)+d^+(x,P)\leq |L|+|E|+1= |P|.   \eqno (23)
$$
 
Finally we consider the paths $H'$ and $Q'$. Now we will compute $d^+(x,H'\cup Q')$ and $d^-(y,H'\cup Q')$.
Denote  $Q'':=Q'\setminus H'$. From the minimality of the path  $H$  it follows that $d^+(x,H')=d^-(y,H')=1$.
This together with the above relations (20)-(23) implies that
$$
d^+(x)+d^-(y)=d^+(x,R)+d^-(y,R)+d^+(x,C[x_{a+1},x_1])+d^-(y,C[x_{a+1},x_1])+d^+(x,C[x_2,x_a])+
$$
$$
d^-(y,C[x_2,x_a])+
d^+(x,P)+d^-(y,P)+
d^+(x,H')+d^-(y,H')+d^+(x,Q'')+d^-(y,Q'')\leq 
$$
 $$
|R|+m-a+1+|P|+2+d^+(x,Q'')+d^-(y,Q'')=
$$ 
$$ 
|R|+|C|+|P|+3+d^+(x,Q'')+d^-(y,Q'')-a. \eqno (24)
$$

Assume that $|H'|\geq 2$, then $ d^+(x,Q'')+d^-(y,Q'')\leq |Q''|$, since otherwise in $\langle B\setminus (P\setminus \{y\})\rangle$ will be an $(x,y)$-path shorter than $H$. The last inequality together with (24) gives
$$
d^+(x)+d^-(y)\leq |R|+|C|+|P|+3+|Q''|-a+|H'|-|H'|=$$ 
$$
n+2-a-|H'|\leq n-2,
$$
since $a\geq 2$ and $|H'|\geq 2$.  This together with (9) implies that
$$d(y)+d(y_1)+d^-(y)+d^+(x)\leq 3n-4,
$$ 
which contradicts  condition $A_0$ since $x, y, y_1$ are $Y$-vertices, $y,y_1$ are nonadjacent and $xy\notin A(D)$.

Now assume that $|H'|=1$, i.e., $|H|=3$, and let $H=xzy$. Now we will compute the sum $d^+(y)+d^-(x)$. It is easy to see that 
$$
d^+(y,R)+d^-(x,R)\leq |R|, \eqno (25)
$$ 
for otherwise in $\langle B\setminus (P\setminus \{y\})\rangle$ there is a $(y,x)$-path which is internally disjoint from $H$.
From the minimality of $C[x_2,x_{a}]$ (i.e., of the gap $a$) it follows that 
$$
d^+(y,C[x_2,x_a])=d^-(x,C[x_2,x_a])=0.   \eqno (26)
$$
 Since the path $H=xzy$ cannot be inserted into $C[x_{a+1},x_1]$ and $x_1x\notin A(D)$, using Lemma 3.2 we obtain 
 $$
d^-(x,C[x_{a+1},x_1])+d^+(y,C[x_{a+1},x_1])\leq m-a+d^-(x,\{x_1\})+d^+(y,\{x_{a+1}\})\leq m-a+1. \eqno (27)
$$
By  minimality of $P$, we have $d^+(y,P)\leq |E|+1$. 
It is not difficult to see that $d^-(x,P)=d^-(x,L) \leq |L|$ (for otherwise using the paths $P_1$, $H$ and the subpaths of the path $P$ we can form a cycle which contains more $Y$-vertices than $C$). The last two inequalities give 
$$
d^+(y,P)+d^-(x,P)\leq |E|+|L|+1=|P|.  \eqno (28)
$$
Note that $H'\cap Q'=\{z\}$.  On the other hand, since $Q$ is minimal, we have $d^+(y,Q')=d^-(x,Q')=1$. Now using this and the above relations  (25)-(28), an easy computation shows that 
$$          
d^+(y)+d^-(x)\leq |R|+|P|+m-a+1+2=|R|+|P|+m+3+|Q'|-|Q'|-a=$$ $$n+2-a-|Q'|\leq n-1,
$$
 since $a\geq 2$ and $|Q'|\geq 1$. This together with (9) contradicts  condition $A_0$, since $x, y, y_1$ are $Y$-vertices, $y, y_1$ are nonadjacent and $yx\notin A(D)$.

 {\it Case 2. In  $\langle B\setminus (P\setminus \{y\})\rangle$ there is no path between the vertices $x$ and $y$. In particular, $x$ and $y$ are nonadjacent.} 

Let now $R:=B\setminus (P\cup \{x\})$. Then it is easy to see that 
$$
 d(x,R)+d(y,R)\leq 2|R|, 
$$
since  in $\langle B\setminus (P\setminus \{y\})\rangle$ there is no path between $x$ and $y$.  Using  Lemmas 3.1 and 3.2 we obtain that
$$
d(x,P)\leq |P|+1 \quad \hbox{and} \quad d(x,C)\leq m,
$$
since the vertex  $x$ cannot be inserted neither into $P$ nor $C$. 
The last three inequalities  together with (4) and (5) imply that 
$$
d(y)+d(x)\leq 2|R|+|P|+m+m-a+|P|+2\leq 2n-a\leq 2n-2.
$$
This together with (9) contradicts Lemma 3.7, since $\{x,y\}$ and $\{y,y_1\}$  are two distinct pairs of nonadjacent vertices of $Y$. The discussion of Case 2 is completed and with it the proof of the theorem is also completed. \fbox \\\\

\section { Concluding remarks }

Observe that the example of the digraph in Remark 1 is not 2-strongly connected and $|Y|=3$. May be true the following conjecture.\\

\noindent\textbf{Conjecture }. {\it Let $D$ be a  digraph of order $n\geq 4$ and let $Y$ be a nonempty subset of vertices of $D$ which satisfies  condition $A_0$.  Then $D$ has a cycle that contains all the vertices of  $Y$ if either (i) or (ii) or (iii) below is satisfied:

(i) $D$ is 2-strongly connected.

(ii) $D$ is $Y$-strongly connected and $|Y|\geq 4$.

(iii) for any ordered pair of distinct vertices $x, y$ of $Y$ there are two internally disjoint paths from $x$ to $y$ in $D$.}\\

C. Thomassen \cite{[18]} (for $n=2k+1$) and the author \cite{[7]} (for $n=2k$) proved the following theorem below. 
Before stating it we need to introduce some additional notations.\\  
Let $m\geq 2$ be any integer.

Let $H(m,m)$ denote the set of digraphs $D$ of order $2m$ with vertex set $A\cup  B$ such that $\langle A\rangle=\langle B\rangle=K^*_m$, there is no arc from $B$ to $A$, $d^+(x,B)\geq 1$ and $d^-(y,A)\geq 1$ for every vertices $x\in A$ and $y\in B$.\\

Let $H(m,m-1,1)$ denote the set of digraphs $D$ of order $2m$ with vertex set $A\cup  B\cup \{a\}$ such that $|A|=|B|+1=m$, $\langle B\rangle \subseteq K^*_m$ (i.e., $\langle B\rangle$ is an arbitrary digraph of order $m-1$) the subdigraph $\langle A\rangle $ has no arc,  $D$ contains all possible arcs betveen $A$ and $B$,  and either $N^-(a)=B$ and $A\subseteq N^+(a)$, or $N^+(a)=B$ and $A\subseteq N^-(a)$.\\

Let $H(2m)$ denote a digraph of order $2m$ with vertex set $A\cup B \cup \{x,y\}$ such that $\langle A\cup \{x\}\rangle=\langle B\cup \{y\}\rangle =K^*_m$, there is no arc between $A$ and $B$, $H(2m)$ also contains all the arcs of the form $ya$, $bx$, for all $a\in A$ and $b\in B$, and the arc $xy$ or both arcs $xy$ and $yx$.\\ 

Denote by $H= H(m,m)\cup H(m,m-1,1)\cup  \{ H(2m)\}$. It is not difficult to show that if a digraph $D$ is in $H$, then its converse digraph also is in $H$ and $D$ is not hamiltonian.\\

Let $D_6$ be a digraph with vertex set $\{x_1,x_2,\ldots ,x_5,x\}$  and arc set 
$$
\{x_ix_{i+1} \,/ 1\leq i\leq 4 \}\cup \{xx_i/\,1\leq i\leq 3\} \cup \{x_1x_5,x_2x_5,x_5x_1,x_5x_4,x_3x_2,x_3x,x_4x_1,x_4x \}.$$ 
By $D'_6$ we denote a digraph obtained from $D_6$ by adding the arc  $x_2x_4$. 

 Note that the digraphs  $D_6$ and $D'_6$ both are 2-strongly connected and are not Hamiltonian. Each  of $D_6$ and $D'_6$ contains a cycle of length 5. Moreover, $d(y)=5$ for any vertex in $D_6$, i.e., $D_6$ is 5-regular. \\ 

\noindent\textbf{Theorem H} ({\it Thomassen \cite{[18]}, Darbinyan \cite{[7]}}). {\it Let  $D$ be a digraph of order $n\geq 5$ with minimum degree at least $n-1$ and with minimum semi-degree at least $n/2-1$. Then $D$ is Hamiltonian unless 

(i) $D$ is isomorphic to $D_5$ or $D_7$ or

$D=[(K_{m}\cup K_{m})+K_1]^*$ or $K^*_{m,m+1}\subseteq D \subseteq [K_{m}+\overline K_{m+1}]^*$, if $n=2m+1$;

(ii) $D\in H\cup  \{ D_6,\, D'_6,\overleftarrow{D_6},\, \overleftarrow {D'_6}\}$, if $n=2m$. 
(The digraphs $D_5$ and $D_7$ are well known and for their definitions,  see, for example \cite{[18]}).} \fbox \\\\

   A question  was put in \cite{[9]}: 
  
Let $D$ be a digraph of order  $n\geq 5$ and let $T\not=\emptyset$ be a subset of $V(D)$. Assume that $D$ is strongly connected (or $D$ is $T$-strongly connected) and every vertex of $T$ has degree at least $n-1$ and has outdegree and indegree at least $n/2-1$. Whether $D$ has a cycle that contains all the vertices of $T$.\\ 

For $n=2m+1$ in \cite{[9]} it was proved:

 {\it If $D$ is strongly connected and contains a cycle of length $n-1$, then $D$ has a cycle containing all the vertices  of $T$ unless some extremal cases.}\\

\end{document}